\documentclass[12pt]{amsart}
\usepackage[all]{xy}

\font\azbuka=wncyr8
\font\itazbuka=wncyi8

\newtheorem{thm}{Theorem}[section]
\newtheorem{prob}[thm]{Problem}

\newcommand{\upannouncement}[1]{[\S\ref{#1} above]}
\newcommand{\dnannouncement}[1]{[\S\ref{#1} below]}
\newcommand{\E}{\exists}
\newcommand{\M}{\mathcal{M}}

\newcommand{\cov}{\mathsf{cov}}

\newcommand{\cof}{\mathsf{cof}}

\newcommand{\spst}{\supseteq}
\renewcommand{\c}{\mathfrak{c}}
\newcommand{\CH}{the Continuum Hypothesis}
\newcommand{\R}{\mathbb{R}}

\newcommand{\cl}[1]{\overline{#1}}

\renewcommand{\b}{\mathfrak{b}}

\newcommand{\bq}{\begin{quote}}
\newcommand{\eq}{\end{quote}}
\renewcommand{\O}{\mathcal{O}}
\newcommand{\B}{\mathcal{B}}
\newcommand{\BG}{\B_\Gamma}

\newcommand{\BO}{\B_\Omega}

\newcommand{\sone}{\mathsf{S}_1}    
\newcommand{\ufin}{\mathsf{U}_{fin}}
    
\newcommand{\seq}[1]{\{#1\}_{n\in\N}}

\newcommand{\cF}{\mathcal{F}}

\newcommand{\NN}{{{}^{\naturals}\naturals}}
\newcommand{\naturals}{{\mathbb N}}
\newcommand{\N}{\naturals}
\newcommand{\sm}{\setminus}
\newcommand{\sbst}{\subseteq}
\newcommand{\by}[2]{\par\hfill\emph{#1}, #2}

\newcommand{\CE}{\textsc{CE}}

\newcommand{\be}{\begin{enumerate}}
\newcommand{\ee}{\end{enumerate}}
\newcommand{\bi}{\begin{itemize}}
\newcommand{\ei}{\end{itemize}}
\renewcommand{\i}{\item}
\newcommand{\SPMBul}{\textbf{$\mathcal{SPM}$ BULLETIN}}

\newcommand{\arx}[1]{\texttt{http://arxiv.org/abs/#1}}

\title[\SPMBul{} \textbf{3} (March 2003)]{%
\SPMBul\\[0.5cm]
Issue number 3: March 2003 \CE{}}

\begin{document}
\maketitle

\tableofcontents

\section{Editor's note}
In this issue we announce a fascinating series of works on
the comparison of various types of convergence of sequences
of functions. Some of these properties are provably
related to some of the properties which were introduced in
the earlier issues of the \SPMBul{}, and many problems remain open.
Section \ref{LBopenprobs} below, written by Lev Bukovsk\'y,
contains a brief survey of some of the major open problems in this area.

This issue gives the first example of the importance
of the transmission of knowledge between the recipients
of this bulletin: One of the announcements implies a solution
to one of the problems posed in an \emph{independent} paper announced
here (see \dnannouncement{edrem}).

The first issues of this bulletin are available online:
\be
\i First issue: \arx{math.GN/0301011}
\i Second issue: \arx{math.GN/0302062}
\ee
We are looking forward to receive more announcements from
other recipients of the bulletin.

\by{Boaz Tsaban}{tsaban@math.huji.ac.il}

\hfill \texttt{http://www.cs.biu.ac.il/\~{}tsaban}

\section{Not distinguishing convergences: Open problems}\label{LBopenprobs}
{
\newcommand{\Zero}[1]{\mbox{\boldmath Z}(#1)}
\newcommand{\wQN}{\hbox{\rm wQN}}
\newcommand{\mQN}{\hbox{\rm mQN}}
\newcommand{\MSS}{\hbox{\rm MSS}}
\newcommand{\QQ}{\mathbb Q}
\newcommand{\RR}{\mathbb R}
We recall some definitions (see e.g.\ \dnannouncement{BuHa}).
An open cover $\mathcal{U}$ is \emph{a $\gamma$-cover} if every
point $x\in X$ is in all but finitely many sets from $\mathcal{U}$.
In accordance with W.\ Hurewicz \cite{Hu} we define \dnannouncement{BuHa}:
\newline
{\bf E}$^*$: {\it for every sequence $\seq{\mathcal{U}_n}$
of open covers of $X$ there exist finite subsets
\mbox{$\mathcal{ V}_n\subseteq\mathcal{ U}_n$} such that
$\seq{\bigcup\mathcal{ V}_n}$ is a cover of $X$}.
\footnote{This is $\ufin(\O,\O)$ in Scheepers' terminology
adopted in this bulletin (see first issue).}
\newline
{\bf E}$^{**}_{\omega}$: {\it for every sequence
$\seq{\mathcal{ U}_n}$ of countable open covers of $X$
there exist finite subsets $\mathcal{ V}_n\subseteq\mathcal{ U}_n$
such that $\seq{\bigcup\mathcal{ V}_n}$ is
a $\gamma$-cover of $X$ or a finite cover of $X$.}
\footnote{This is $\ufin(\O,\Gamma)$ in Scheepers' terminology,
if $\O$ is restricted to \emph{countable} open covers of $X$.
This restriction can make a difference when $X$ is not Lindel\"of.}

A countably compact non-compact topological spa\-ce has property
E${}^{**}_{\omega}$ and has not  property E${}^*$. There are examples
of such spaces {\rm (e.g.\ \cite{En}, pp.\ 261--262)}, however none
of them is perfectly normal.
The existence of a perfectly normal countably compact non-compact
space neither can be proved nor can be refuted in {\bf ZFC}
{\rm (see e.g.\ \cite{Va})}.
\begin{prob}[\dnannouncement{BuHa}]\label{1}
Find in ZFC a perfectly normal
{\rm E}${}^{**}_{\omega}$-space which does not possess
property~{\rm E}${}^*$.
\end{prob}
\par
We say that a sequence $\seq{f_n}$
converges \emph{quasi-normally to a function $f$ on $X$},
(see e.g.\ \cite{Ba}, in \cite{CL} as equally convergent)
if there is a sequence of positive reals $\seq{\varepsilon_n}$
(\emph{a control}) converging to $0$ such that
\begin{equation}\label{quasi}
(\forall x\in X)(\exists n_0)(\forall n\geq n_0)\,
\vert f_n(x)-f(x)\vert <\varepsilon_n.
\end{equation}
Similarly, the series $\sum_{n=0}^{\infty}f_n$ converges
\emph{pseudo-normally on $X$} if there is a control sequence
$\seq{\varepsilon_n}$ such that
$\sum_{n=0}^{\infty}\varepsilon_n<\infty$ and (\ref{quasi})
holds true (with $f=0$).
\par
A topological space $X$ is said to be \emph{a wQN-space},
see \cite{BRR1}, if
from every sequence of continuous functions converging to 0
on $X$ one can choose a quasi-normal\-ly convergent subsequence.
A topological space $X$ is said to
be \emph{a} $\Sigma\Sigma^*$\emph{-space}, see \dnannouncement{BRR2}, if
for every sequence
$\seq{f_n}$ of real functions with non-negative values such that
$\sum_{n=0}^{\infty}f_n(x)<\infty$ for every $x\in X$ the series
converges also pseudo-normally.
Finally a topological space is \emph{a $\overline{\rm QN}$-space},
see \dnannouncement{BRR2}, if every sequence of real functions converging
pointwise to a function on $X$ (not necessarily continuous) converges
to this function quasi-normally.
\par
In \dnannouncement{BRR2} the authors show that
\[\Sigma\Sigma^*\to{\overline{\rm QN}}\mbox{\ for perfectly normal
space}.\]
Usually passing from properties of sequences of real-valued continuous
functions to properties of open coverings we need to assume that the
considered topological space is perfectly normal. However, both notions
$\Sigma\Sigma^*$-space and ${\overline{\rm QN}}$-space do not use
a notion of an open covering in their definitions. Therefore
we suppose that
\begin{prob}\label{2}
$\Sigma\Sigma^*\to{\overline{\rm QN}}$ for
arbitrary topological space.
\end{prob}
For a topological space $X$ and a subset $A\subset X$ we denote
\begin{eqnarray*}
s_0(A)&=&A,\quad s_{\xi}(A)=\{\lim_{n\to\infty} x_n:x_n\in
\bigcup_{\eta<\xi}s_{\eta}(A)\mbox{\ for each\ }n\in\N\},\\
\sigma(A)&=&\min\{\xi:s_{\xi}(A)=s_{\xi+1}(A)\},\quad
\Sigma(X)=\sup\{\sigma(A) : A\subseteq X\},
\end{eqnarray*}
The fundamental result in this area is David Fremlin's
\begin{thm}[\cite{Fr1}] $\Sigma(C_p(X))=0,\,1,\,\omega_1$.
\end{thm}
The theorem suggests to define: a topological space $X$ is said to be
\emph{an $s_1$-space} if $\Sigma(C_p(X))=1$.
\par
In \cite{Sc1} the author introduces the \emph{sequence selection
property}, shortly \emph{SSP} of a topological space $X$:
if $\lim_{i\to{\infty}}f_{n,i}(x)=0$ for $x\in X$,
$n\in\N$, then there are~$i_n$ such that
$\lim_{n\to{\infty}}f_{n,i_n}(x)=0$ for $x\in X$.
Actually SSP is equivalent to $\alpha_2$
property of $C_p(X)$ introduced by A.\ V.\ Ar\-chan\-gel\-skij~
\cite{Ar}.
\begin{thm}{\rm (\cite{wqn}, implicitly in \cite{Fr1})}
\ {\rm SSP} = $s_1$-space.
\end{thm}
\begin{thm}{\rm (\cite{wqn})}\ {\rm SSP} $\to$ {\rm wQN}.
\end{thm}
Recently D.\ Fremlin  proved
\begin{thm}{\rm (\cite{Fr2})}{\rm\ wQN$\to$SPP}.
\end{thm}
A topological space $X$ is said to be
\emph{a ${\rm S}_1(\Gamma,\Gamma)$-space} if for every sequence
$\seq{\mathcal{U}_n}$ of $\gamma$-covers of $X$ there exists
a~$\gamma$-cover $\seq{U_n}$ such that $U_n\in\mathcal{U}_n$
for every $n\in\N$. In \cite{wqn} the author shows that
$\sone(\Gamma,\Gamma)\to {\rm s}_1{\rm -space}$
and conjectured that
\begin{prob}\label{3}
Every perfectly normal wQN-space {\rm (= s${}_1$-space)} has
property $\sone(\Gamma,\Gamma)$.
\end{prob}

\by{Lev Bukovsk\'y}{bukovsky@kosice.upjs.sk}
}

\section{Research announcements}

\subsection{Spaces not distinguishing convergences of real-valued functions}\label{BRR2}
In \cite{BRR1} we have introduced the notion of
a wQN-space as a~space in which for every
sequence of continuous functions
pointwisely converging to $0$
there is a subsequence
quasi-normally converging to $0$.
In the present paper we continue this investigation
and generalize some concepts touched there.
The content is a~variety of notions and
relationships among them.
The result is another scale in the investigation
of smallness and the question is how this scale
fits with other known scales and whether all
relations in it are proper.

\noindent
The paper appeared in \emph{Topology and its Applications} \textbf{112} (2001), 13--40.

\by{Lev Bukovsk\'y}{bukovsky@kosice.upjs.sk}
\by{Ireneusz Rec\l{}aw}{reclaw@ksinet.univ.gda.pl}
\by{Miroslav Repick\'y}{repicky@kosice.upjs.sk}

\subsection{Hurewicz Properties, not Distinguishing Convergence
Properties and Sequence Selection Properties}
We shall compare several properties of a topological space
related to the behavior of open coverings and/or the
behavior of sequences of continuous real-valued functions defined
on the space. We shall show that there are closed relationships
between them and several of them are mutually equivalent.
\by{Lev Bukovsk\'y}{bukovsky@kosice.upjs.sk}

\subsection{On Hurewicz Properties}\label{BuHa}
We investigate  Hurewicz properties introduced in \cite{Me} and
\cite{Hu} and later introduced related properties of topological
spaces. The main result says that for perfectly normal spaces
the property mQN introduced in \upannouncement{BRR2} is equivalent
to Hurewicz property E${}^{**}_{\omega}$. As corollaries we obtain
solution of several open problems stated in \upannouncement{BRR2}.
A complete overview of relationships between the considered
properties is presented.
\by{Lev Bukovsk\'y}{bukovsky@kosice.upjs.sk}
\by{Jozef Hale\v s}{hales@science.upjs.sk}

\subsection{Uncountable $\Sigma\Sigma^*$ subset of $\R$}
If CH holds then there exists an uncountable $X\subset[0,1]$
which belongs to $\Sigma\Sigma^*$.
\by{Lev Bukovsk\'y}{bukovsky@kosice.upjs.sk}
\by{Krzysztof Ciesielski}{K\_Cies@math.wvu.edu}

\subsection{Spaces not distinguishing convergences}
In the present paper we introduce a~convergence
condition ($\Sigma'$) and continue the study of
``not distinguish'' for various kinds of
convergence of sequences of real functions on
a~topological space started in \cite{BRR1}
and \upannouncement{BRR2}.
We compute cardinal invariants associated with
introduced properties of spaces.
\by{Miroslav Repick\'y}{repicky@kosice.upjs.sk}

\subsection{A Nonhereditary Borel-cover $\gamma$-set}\label{NonHered}
In this paper we answer some of the questions raised
by Bartoszynski and Tsaban \cite{ideals} concerning hereditary
properties of sets defined by certain Borel covering properties:\\
\textbf{Theorem.} Suppose there is a Borel-cover $\gamma$-set
\footnote{That is, an element of $\sone(\BO,\BG)$.}
of size the continuum. Then there is a Borel-cover $\gamma$-set $X$ and subset $Y$ of $X$ which is
not even an open-cover $\gamma$-set.  (In fact there is an
of open $\omega$-cover of $Y$ with no $\tau$-subcover.)\\
It is also shown that CH implies that there exists a Borel-cover $\gamma$-set of size
$\omega_1$.
\by{Arnold W.\ Miller}{miller@math.wisc.edu}

\subsection{Editor's remark: A $\gamma$- and $\sigma$-set need not be hereditary}\label{edrem}
In Problem 7.9 of the announced paper \upannouncement{BRR2}
it is asked whether every $\gamma$-set of reals
which is also a $\sigma$-set is a \emph{hereditary} $\gamma$-set.
By \upannouncement{NonHered}, assuming CH there exists
an element of $\sone(\BO,\BG)$ with a subset which is not a $\gamma$-set.
Clearly $\sone(\BO,\BG)$ implies $\sone(\Omega,\Gamma)$ (= $\gamma$-set),
as well as $\sone(\BG,\BG)$. In \cite{CBC} it is proved that
every set satisfying $\sone(\BG,\BG)$ is a $\sigma$-set.
This answers the problem negatively.
\by{Boaz Tsaban}{tsaban@math.huji.ac.il}

\subsection{The minimal cardinality where the Reznichenko property fails}
According to Reznichenko,
a topological space $X$ has the weak Fr\'echet-Urysohn
property if for each subset $A$ of $X$ and each element
$x$ in $\cl{A}\sm A$, there exists a countably infinite
pairwise disjoint collection $\cF$ of finite subsets of $A$
such that for each neighborhood $U$ of $x$, $U\cap F\neq\emptyset$ for
all but finitely many $F\in\cF$.
In \cite{FunRez}, Ko\v{c}inac and Scheepers conjecture:
\begin{quote}
The minimal cardinality of a set $X$ of real numbers
such that $C_p(X)$ does not have the weak Fr\'echet-Urysohn property
is equal to $\b$.
\end{quote}
($\b$ is the minimal cardinality of an unbounded family in the
Baire space $\NN$).
We prove the Ko\v{c}inac-Scheepers conjecture by showing that
if $C_p(X)$ has the Reznichenko property, then a continuous image of
$X$ cannot be a subbase for a non-feeble filter on $\N$.
\by{Boaz Tsaban}{tsaban@math.huji.ac.il}

\section{Other announcements}

\subsection{The Twentyseventh Summer Symposium in Real  Analysis}
June 23--29, 2003.

During June 23-29, 2003, the Mathematical Institute of Silesian
University at Opava will host the Summer Symposium in Real
Analysis XXVII. The nature of current work in real analysis is
driven by the exchange of ideas generated by real analysts rooted
in one subdiscipline of real analysis but with wide ranging
interests. This Symposium will highlight lectures by both leading
experts and energetic new researchers. Specifically, Summer
Symposium XXVII will emphasize recent important work in harmonic
analysis, integration theory and a solution of the celebrated
Gradient Problem as well as some of the achievements of younger
mathematicians in real analysis. In addition, we will provide a
vibrant forum for the discussion of research problems, and allot
prime speaking time to recent doctoral recipients.

The principal speakers have been invited and at the time of
this submission, all have tentatively accepted our invitation.
\bi
\i Jaroslav Kurzweil (Mathematical Institute of the Academy of
Sciences, Prague)
\i Zoltan Buczolich (Eotvos Lorand University, Budapest)
\i Alexander Olevskii (Tel Aviv University, Israel)
\ei
Michigan State University Press will publish the proceedings of Symposium XXVI as a
separate volume of the Real Analysis Exchange.
Electronic registration for the conference can be found at:
\bq
\texttt{http://www.math.slu.cz/RealAnalysis/}
\eq
\by{Petra Sindelarova}{Petra.Sindelarova@math.slu.cz}

\subsection{BEST 2003 (update)}
The organizers of the BEST 2003 conference
have informed us with the good news that Arnold Miller will attend
this conference and give an invited lecture.
For more details see \cite{BEST03}.

\section{Problem of the month}
The following problem, which is a variant of Problem 3 in \cite{coc2},
appears as Problem 1 in \upannouncement{BuHa} and in \upannouncement{LBopenprobs}.
\begin{prob}
Does there exist (in ZFC) a set $X\sbst\R$ which has the Menger property
$\ufin(\O,\O)$ but not the Hurewicz property $\ufin(\O,\Gamma)$?
\end{prob}
Assuming \CH{}, one can construct
a \emph{Luzin set} $L\sbst\R$ of size continuum $\c$,
that is, such that for each meager (=first category)
set $M$, $L\cap M$ is countable. Such a set $L$ is
\emph{concentrated} on each of its countable dense subsets $D$
(that is, for each open set $U\spst D$, $L\sm U$ is countable),
and therefore has Rothberger's property $\sone(\O,\O)$, which
implies Menger's property $\ufin(\O,\O)$ (see, e.g., \cite{MilSpec}).
On the other hand, in \cite{coc2} it is proved that every set with
the Hurewicz property is (perfectly) meager.
Thus, assuming \CH{}, the answer to the above problem is negative.
But \CH{} is not necessary to get a negative answer:
Let $\M$ denote the collection of meager sets of reals, and
write
\begin{eqnarray*}
\cov(\M) & = & \min\{|\cF| : \cF\sbst\M \mbox{ and }\cup\cF=\R\}\\
\cof(\M) & = & \min\{|\cF| : (\forall M\in\M)(\E F\in\cF)\ M\sbst F\}
\end{eqnarray*}
In \cite{CBC} it is shown that it is enough to assume that
$\cov(\M)=\cof(\M)$ (this hypothesis is strictly weaker than \CH{} \cite{jubar})
for the above arguments to work (with some necessary modifications).
The Problem of the Month asks whether the assumption
$\cov(\M)=\cof(\M)$ can be completely removed.

The papers \cite{coc2}, \cite{wqn}, and \cite{ideals}
deal with constructions in ZFC of sets of reals with the Hurewicz property,
and seem to be relevant to the problem.
\by{Boaz Tsaban}{tsaban@math.huji.ac.il}

\end{document}